\documentclass[12pt]{article}  

\usepackage{amsmath}
\usepackage{amssymb}
\usepackage{cite}
\usepackage{paralist}
\usepackage[usenames]{color}
\usepackage{graphicx}          
\usepackage[dvips]{epsfig}    
\newtheorem{definition}{Definition}
\newtheorem{theorem}{Theorem}
\newtheorem{corollary}{Corollary}
\newtheorem{remark}{Remark}
\newtheorem{problem}{Problem}

\newtheorem{assumption}{Assumption}

\title{Partial Stability Concept in Extremum Seeking Problems
\thanks{
This work is  supported in part by the German Research Foundation (GR 5293/1-1).\newline
$^{1}$Institute of Mathematics, Julius Maximilian University of Würzburg,   Germany
        {\tt\small viktoriia.grushkovska@mathematik.uni-wuerzburg.de}
        \newline
$^{2}$Max Planck Institute for Dynamics of Complex Technical Systems, Magdeburg, Germany
{\tt\small zuyev@mpi-magdeburg.mpg.de}\newline
$^{3}$Institute of Applied Mathematics and Mechanics, National Academy of Sciences of Ukraine
}
}

\author{Victoria Grushkovskaya$^{1,3}$ and Alexander Zuyev$^{2,3}$}
\date{}

\begin{document}

\maketitle
\thispagestyle{empty}

\begin{abstract}
The paper deals with the extremum seeking problem for a class of cost functions depending only on a part of state variables  of a control system.
This problem is related to the concept of partial asymptotic stability and analyzed by Lyapunov's direct method and averaging schemes.
Sufficient conditions for the practical partial stability of a system with oscillating inputs are derived with the use of Lie bracket approximation techniques.
These conditions are exploited to describe a broad class of extremum-seeking controllers ensuring the partial stability of the set of minima of a cost function.
The obtained theoretical results are illustrated by the Brockett integrator and rotating rigid body.
\end{abstract}

\section{Introduction}

Extremum seeking has become an important branch of modern control theory because of challenging theoretical features and various practical applications. The goal of extremum seeking control is to optimize the steady-state performance of a control system using  the output measurements. The main motivation behind this problem statement is to reduce the amount of information needed for the control design. In particular, an optimal operating point as well as analytical expression of the output (cost) function are assumed to be unknown. During the past couple of decades, several important approaches for the extremum seeking control design have been developed (see, e.g.,~\cite{Kr00,Krst03,Guay03,review,Nes10,Durr13a,Guay15,Har16,Durr17,SD17,SK17,GMZME18,GZE18}). The above approaches  assume that the cost function depends essentially on all state variables, and/or that the system admits an asymptotically stable steady-state. However, these assumptions can be redundant for various applied problems, for which it is important (or even only possible) to optimize the system with respect to a prescribed part of state variables, and consequently to stabilize the system only with respect to these variables. In particular, such problems arise if the cost function depends on a part of system variables, if only  partial output measurements are available for control design, or if the partial stabilization is sufficient for correct system operation. As a simple example, one can imagine the problem of tracking a planar target by a multi-DOF robot (see, e.g., \cite{Coh09,Mat11,Khong14,Mand15})
\\
The goal of this paper is to introduce the problems of partial extremum seeking, in which the goal is to optimize the system performance with respect to a part of state variables only.
 Such problem statement allows to consider a broader class of systems and applications.
 The contribution of this paper is twofold. First, we generalize the Lie bracket approximation approach (see, e.g.,~\cite{MA00,Durr13a,Durr17}) and  techniques introduced in~\cite{GZE18} to input-affine systems whose Lie bracket system has a partially asymptotically stable manifold. To solve the problem under consideration, we  attract methods of  partial  stability theory, which dates back to Lyapunov and has been developed in the works of~\cite{Mal,RO87,Vo12,Z2000,Zu03,Kov09,GZ15} and others (see \cite{Vo05} for a review).
 Second, we consider a class of extremum seeking problems, in which the system has to be optimized with respect to a prescribed part of variables. Up to our best knowledge, such problem statement has not been considered before.
\\
 The rest of the paper is organized as follows. Section~1.1 contains some  notations and definitions which will be used throughout the paper. In Section~2.1, we extend the Lie bracket approximation approach assuming that the corresponding Lie bracket system is partially asymptotically stable, and derive conditions  for practical partial asymptotic stability. These results are applied to extremum seeking problems in Section~2.2. In Section~3, we consider several examples  illustrating the proposed approach and some possible extensions.

\subsection{Notations and definitions}
Consider the system
\begin{equation}\label{sys}
  \dot x =f^\varepsilon(t,x),\;x\in\mathbb R^n,\,t\ge 0,
\end{equation}
where $f:\mathbb R^+\times \mathbb R^n \to\mathbb R^n$, and $\varepsilon>0$ is a parameter.   We will split the components of the state vector $x$ as $x=(y^\top,z^\top)^\top\in\mathbb R^n$ with $y\in\mathbb R^{n_1}$, $z\in\mathbb R^{n_2}$, $n_1+n_2=n$.
 With a slight abuse of notations, the column $x$ will be also denoted as $x=(y,z)$.
Throughout the text,
 $B_\delta(x^*)$ and $\overline{B_\delta(x^*)}=B_\delta(x^*)\cup \partial B_\delta(x^*)$ denote the $\delta$-neighborhood of an $x^*{\in} \mathbb R^n$ and its closure, respectively.  Notation $\varphi\in \mathcal K$ means that a function $\varphi$ belongs to the class $\mathcal K$, i.e. $\varphi:\mathbb R^+\to\mathbb R^+$ is a continuous strictly increasing function, $\varphi(0)=0$.  For  $f,g:\mathbb R^n\to\mathbb R^n $, $x\in\mathbb R^n$, we denote the directional derivative as
 $ L_gf(x)=\lim\limits_{s\to0}\tfrac{f(x+sg(x))-f(x)}{s}$, and  $[f,g](x)= L_fg(x)- L_gf(x)$ is the Lie bracket.
We will use the following definition, which extends the notion of \emph{partial asymptotic stability}~(\cite{RO87,Vo12,Zu03}) to systems with parameters of the form~\eqref{sys}.
\begin{definition}
For $y^*\in\mathbb R^{n_1}$, the set $D^*=\{x=(y,z)\in\mathbb R^n: y=y^*\}$ is \emph{practically uniformly $y$-asymptotically stable for system~\eqref{sys}, }if it is:
\\
  $-$ \emph{practically  uniformly $y$-stable for system~\eqref{sys}, }i.e., for every $\rho>0$, there exist $\delta>0$, $\bar\varepsilon>0$ such that the following property holds for all $t_0\ge 0$, $z(t_0)\in \mathbb R^{n_2}$,  $\varepsilon\in(0,\bar\varepsilon)$:
      $$
      \text{if }y(t_0)\in B_\delta(y^*)\text{ then }y(t)\in B_\rho(y^*)\text{ for all }t\in[t_0,\infty);
      $$
  $-$ \emph{practically uniformly $y$-attractive for system~\eqref{sys}, }i.e., for some $\delta{>}0$, for every $\rho{>}0$, there are $t_1{\ge} 0$,  $\bar\varepsilon{>}0$ such that the following property holds for all  $t_0{\ge }0$, $z(t_0)\in \mathbb R^{n_2}$,  $\varepsilon{\in}(0,\bar\varepsilon)$:
      $$
      \text{if }y(t_0)\in B_\delta(y^*)\text{ then }y(t)\in B_\rho(y^*)\text{ for all }t\in[t_0+t_1,\infty).
      $$
If the attractivity property holds for any $\hat\delta{>}0$, then~$y^*$ is called to be \emph{semi-globally practically uniformly $y$-asymptotically stable for system~\eqref{sys}. }
For systems independent of $\varepsilon$, we omit the terms ``practically'' and ``semi''.\\
In case $n_1=n$, $n_2=0$, the above definition coincides with a well-known definition of practical asymptotic stability~(\cite{MA00,Durr13a}).
Up to our best knowledge, the proposed definition of \emph{practical} partial stability is introduced here for the first time.
\end{definition}
\section{Main results}
\subsection{Lie bracket approximation \& partial  stability}
In this section, we extend the Lie bracket approximation approach to partially asymptotically stable systems.
Namely, we consider the system
\begin{equation}\label{sysA}
\begin{aligned}
\dot x = f_0(x)+\sum_{i=1}^mf_i(x)u_i,\;x\in\mathbb R^n,\,
\end{aligned}
\end{equation}
where $u_i{=}\tfrac{1}{\sqrt\varepsilon}w_i\Big(\tfrac{t}{\varepsilon}\Big)$, $w_i\Big(\tfrac{t}{\varepsilon}\Big)$ are $\varepsilon$-periodic continuous functions with some $\varepsilon{>}0$, and $\int_0^\varepsilon w_i\Big(\tfrac{t}{\varepsilon}\Big)dt{=}0$. We assume that there exists a $W{>}0$ such that
$\max\limits_{1\le i\le m,0\le t\le\varepsilon}w_i\Big(\tfrac{t}{\varepsilon}\Big){\le} W\text{ for each } \varepsilon{>}0.$
Consider also the so-called Lie bracket system
 \begin{equation}\label{sysB}
\begin{aligned}
\dot{\bar x} = f_0(\bar x)+\sum_{i<j,i,j=1}^m[f_i, f_j](\bar x)\nu_{ij},\;\bar x\in\mathbb R^n,
\end{aligned}
\end{equation}
where  $\nu_{ij}=\tfrac{1}{\varepsilon^2}\int_0^\varepsilon\int_0^\tau w_j\Big(\tfrac{\tau}{\varepsilon}\Big)w_i\Big(\tfrac{s}{\varepsilon}\Big)dsd\tau$. Denote $x=(y,z)$, $\bar x=(\bar y,\bar z)$, $y,\bar y\in\mathbb R^{n_1}$, $z,\bar z\in\mathbb R^{n_2}$, $n_1+n_2=n$.
\begin{assumption}
Let $D_1\subseteq\mathbb R^{n_1}$ and $D_2\subseteq\mathbb R^{n_2}$ be domains,  and let $y^*\in D_1$, $D=\{(y,z)\in\mathbb R^n:y\in D_1,z\in D_2\},$ $D^*=\{x=(y,z)\in D: y=y^*\}.$
We suppose that:
\begin{itemize}
\item[A1.1)]     $f_0,f_1,\dots,f_m\in C^2(D\setminus D^*;\mathbb R^n)$;
\item [A1.2)] for any compact $\widetilde D_1\subset D_1$, the functions  $f_i$, $L_{f_j}f_i$, $L_{f_l} L_{f_j}f_i\in C(D;\mathbb R^n)$ are bounded for all $y\in \widetilde  D_1$, $z\in D_2$,   $i,j,l\in\{0,\dots,m\}$;
\item[A1.3)] if $x(t){\in }D$, $t{\in} I{=}[t_0,t_1)$ is a solution of~\eqref{sysA} s.t. $\inf\limits_{t\in I}{\rm dist}(y(t),\partial D_1){>}0$ then $\inf\limits_{t\in I}{\rm dist}(z(t),\partial D_2){>}0$.
\end{itemize}
\end{assumption}
Here ${\rm dist}(\xi,X)$ denotes the Euclidian distance between a point $\xi{\in}\mathbb R^{n_k}$ and a set $X{\subset}\mathbb R^{n_k}$. If both $D_2{\subset} \mathbb R^{n_2}$ and $z(t){\in} D_2$, $t{\in} I$, are unbounded, we will follow the convention that $\inf_{t\in I}{\rm dist}(z(t),\partial D_2)=0$.  Note that A1.3) is a reformulation of the standard $z$-extendability assumption in partial stability theory (see, e.g.,~\cite{RO87}). For the case $D_2{=}\mathbb R^{n_2}$, this assumption means that $z(t)$ cannot escape to infinity in finite time whenever $y(t)$ remains bounded. The above assumption is usually satisfied in well-posed practical problems without blow-up of solutions.\\
 The first main result of the paper is as follows.
\begin{theorem}~\label{thm_gen_practical}
 {Let $D_1\subseteq\mathbb R^{n_1}$, $D_2\subseteq\mathbb R^{n_2}$ be such that Assumption~1 is satisfied, $y^*\in D_1$, and  let there exist a function $V(x)\in C^2(D)$ such that the following conditions hold for all $x=(y,z)\in D$:
  \begin{itemize}
    \item[1.1)] $\alpha_{1}(\|y-y^*\|)\le V(x)\le\alpha_{2}(\|y-y^*\|)$,
    \item[1.2)] $ L_{\bar f}V(x)\le -\alpha_3(\|y-y^*\|)$.
    \end{itemize}
Here $\bar f(x)=f_0(x)+\sum_{i<j}[f_i, f_j](x)\nu_{ij}$  is the right-hand side of system~\eqref{sysB}, and  $\alpha_1,\,\alpha_2,\,\alpha_3\in \mathcal K$.\\
   Then $D^*=\{x=(y,z)\in D: y=y^*\}$ is   practically $y$-asymptotically stable for~\eqref{sysA} with the initial conditions from the set  $D_0=\{(y,z)\in\mathbb R^n:\|y-y^*\|\le \delta,z\in D_2\}$, where $\delta\in\Big(0,\alpha_2^{-1}\big(\alpha_1({\rm dist}(y^*,\partial D_1))\big)\Big)$.}
\end{theorem}
The proof of Theorem~\ref{thm_gen_practical} is in Appendix~A. Note that the assumptions of Theorem~\ref{thm_gen_practical} are more general than those used in~\cite{GZE18}, so that the proof of this result extends the approaches of~\cite{GZE18} to a broader class of systems. \\
The next results follow from the proof of Theorem~1.
\begin{corollary}
{If the conditions of Theorem~\ref{thm_gen_practical} hold with~1.1) replaced by $\alpha_{1}(\|y-y^*\|)\le V(x)\le\tilde\alpha_{2}(\|x-x^*\|)$, $x^*=(y^*,z^*)$, $z^*\in D_2$, where $\tilde\alpha_2\in \mathcal K$,  then the set $D^*$ is   practically $y$-attractive in $D_0$ for system~\eqref{sysA} provided that there exist  $\delta>0$, $c_{\delta}\in\Big(0,\alpha_1\big({\rm dist}(y^*,\partial  D_1)\big)\Big)$ such that $\alpha_2(\|x-x^*\|)\le c_{\delta}$ for all  $x\in D_0$.}
\end{corollary}
\begin{corollary}
{If the conditions of Theorem~\ref{thm_gen_practical} hold with the function $V$ depending on the $y$-variable only,  then the assertion of Theorem~\ref{thm_gen_practical} holds even if the $z$-components of the
 functions from~A1.2) are unbounded.}
\end{corollary}
\begin{remark}
Under some additional assumptions on the function $V$ and the vector fields of system~\eqref{sysA}, it is possible to state \emph{classical} (instead of practical) asymptotical stability conditions and to describe the decay rate of solutions of system~\eqref{sysA}, as it was done in~(\cite{GZE18}) by extending the techniques of (\cite{GZ14,ZG17}). We leave these studies for future work.
\end{remark}

\subsection{Partial stabilization of control-affine \\ extremum seeking systems}
In this section, we apply the proposed results to extremum seeking problems in which the goal is to optimize the system performance with respect to certain part of variables.
Namely, we assume that the set of minima of a cost function $J:\mathbb R^n \to \mathbb R$ is a hyperplane of the form
${\rm argmin}\, J = \{x=(y,z):y=y^*\}$, where the value of $y^*\in \mathbb R^{n_1}$ is a priori unknown for the control design.
Thus we arrive to the following problem statement.
\begin{problem}\label{problem-statement}
 Given a cost function $J{:}\,\mathbb R^{n}\to\mathbb R$ such that
 $${\rm argmin}\, J {=} \{x{=}(y,z){\in}\mathbb R^n:y{=}y^*\}\;\; \text{with some}\; y^*{\in} \mathbb R^{n_1}.$$
The goal is to {construct a control $u=\tilde u(t,J(x))$ such that the set  ${\rm argmin}\, J$ is   practically $y$-asymptotically stable for~\eqref{sysA}.}
\end{problem}
Such kind of problems appears, for example, if the cost $J$ depends on the $y$-variables only,
or if $J$ can be represented as $J(x)=J^*(y-y^*)\phi(z)$, where $J^*(\eta)$ is a positive definite function, and $\phi(z)>0$ for all $z\in\mathbb R^{n_2}$.
The above task is relevant to the output stabilization problem, if the stabilization with respect to all variables is not possible (or not required for control purposes),
and to synchronization problems, where the goal $y=y^*$ describes synchronous motion of a multi-agent system (e.g., system of pendulums) while the $z$-variables stand for redundant degrees of freedom.
Let us define the controls $u_i$ as
\begin{equation}\label{cont}
\begin{aligned}
  u_i=\frac{1}{\sqrt{\varepsilon}}\Big(g_i(J(x))w_i\Big(\frac{t}{\varepsilon}\Big){+}g_{i+m}(J(x))w_{i+m}\Big(\frac{t}{\varepsilon}\Big)\Big),
\end{aligned}
\end{equation}
where $\varepsilon>0$, $w_i,w_{i+n}$ satisfy the assumptions of section~2.1 and are such that $\nu_{ij}=0$ whenever $j\ne i+m$, $\nu_{ii+m}=1$, and the functions $g_i,g_{i+m}$ satisfy the  relation
\begin{equation}\label{class}
g_{i+m}(z)=-\gamma_i g_{i}(z)\int{\frac{dz}{g_{i}(z)^2}}, \gamma_i>0, \ i=\overline{1,m}.
\end{equation}
\begin{theorem}~\label{thm_es}
{Let $D_1\subseteq\mathbb R^{n_1}$, $D_2\subseteq\mathbb R^{n_2}$ be convex domains such that Assumption~1 is satisfied, $y^*\in D_1$, and let the function $V(x)=J(x)-J(y^*,z)$ satisfies the conditions of Theorem~\ref{thm_gen_practical} with
 $$\bar f(x)=f_0(x)-\sum_{i=1}^m\gamma_if_i(x)f_i^\top(x)\nabla J(x).$$
    Then the set $D^*=\{x=(y,z)\in D: y=y^*\}$ is practically $y$-asymptotically stable in $D_0$ for system~\eqref{sysA} with the controls $u_i$ given by~\eqref{cont}--\eqref{class}.
}
\end{theorem}
{\bf Proof.}
 Straightforward calculations show that the Lie bracket system for~\eqref{sysA} with the controls $u_i$ given by~\eqref{cont}--\eqref{class} has the form
  $$
  \dot{\bar x}=f_0(\bar x)-\sum_{i=1}^m\gamma_if_i(\bar x)f_i^\top(\bar x)\nabla J(\bar x).
  $$
Then the conditions of Theorem~\ref{thm_gen_practical} are satisfied.
$\square$

The  assumptions on the cost function $J$ required in Theorem~\ref{thm_es} are  common in extremum seeking studies for ensuring the stability with respect to all variables (cf.~\cite{Tan06,Gu18}).
They can be relaxed for certain classes of systems, as in the next result.
\begin{theorem}~\label{thm_es2}
{Let a control system be of the form
  \begin{equation}\label{sys_es_1}
\begin{aligned}
   & \dot y=\sum_{i=1}^{n_1}\tilde f_i(x)u_i,\quad \dot z= h(x,u),
\end{aligned}
  \end{equation}
  where the vector fields $\tilde f_i:\mathbb R^{n}\to\mathbb R^n$ and $h:\mathbb R^n\times\mathbb R^{n_1}\to\mathbb R^{n_2}$ satisfy A1.1)--A1.2). Assume that the vector fields $ \tilde f_i(x)=(\tilde f_{i1}(x)   \dots   \tilde f_{in_1}(x))^\top$, $i=1,2,\dots,n_1$, are linearly independent at each $x\in D$, and the cost function $J=J(y):D_1\subset\mathbb R^{n_1}\to\mathbb R$ satisfies the inequalities
$$\alpha_{1}(\|y-y^*\|)\le J(y)-J(y^*)\le\alpha_{2}(\|y-y^*\|),$$
 $$\|\nabla J(y)\|\le -\alpha_3(\|y-y^*\|)$$ with some $ \alpha_1,\alpha_2,\alpha_3\in \mathcal K$.\\
Then the set $D^*$ is   practically $y$-asymptotically stable  for system~\eqref{sys_es_1} with the controls $u_i$ given by~\eqref{cont}--\eqref{class}.}
\end{theorem}
\begin{figure*}[ht]
\begin{minipage}{0.32\linewidth}
\begin{center}
\includegraphics[width=1\linewidth]{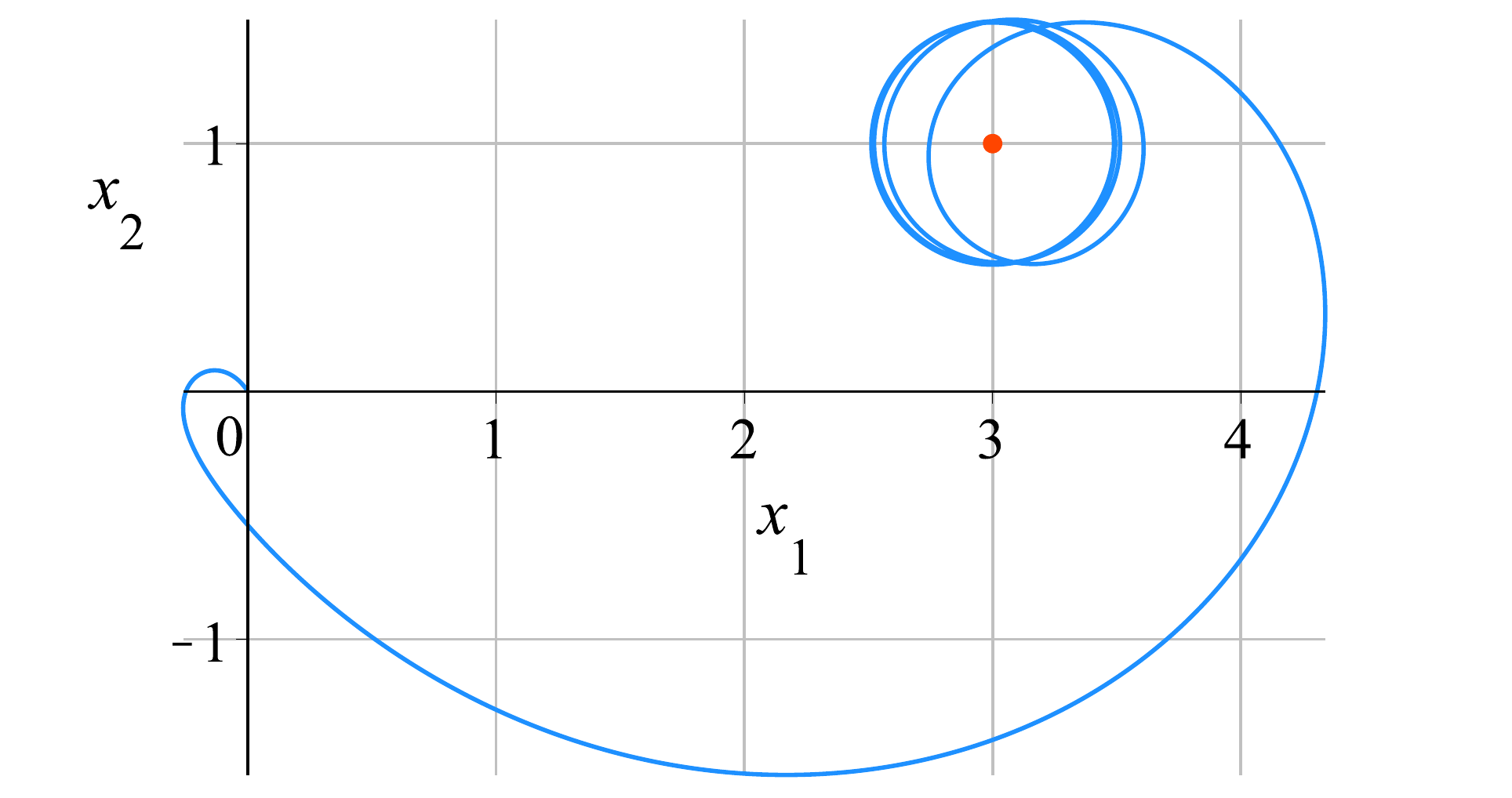}\\
\includegraphics[width=1\linewidth]{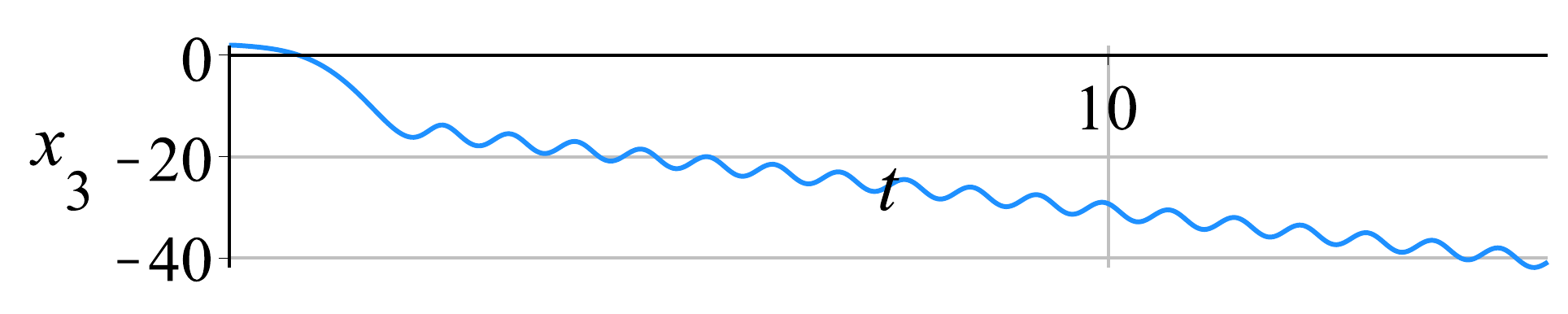}\\
a)
\end{center}
\end{minipage}
\begin{minipage}{0.32\linewidth}
\begin{center}
\includegraphics[width=1\linewidth]{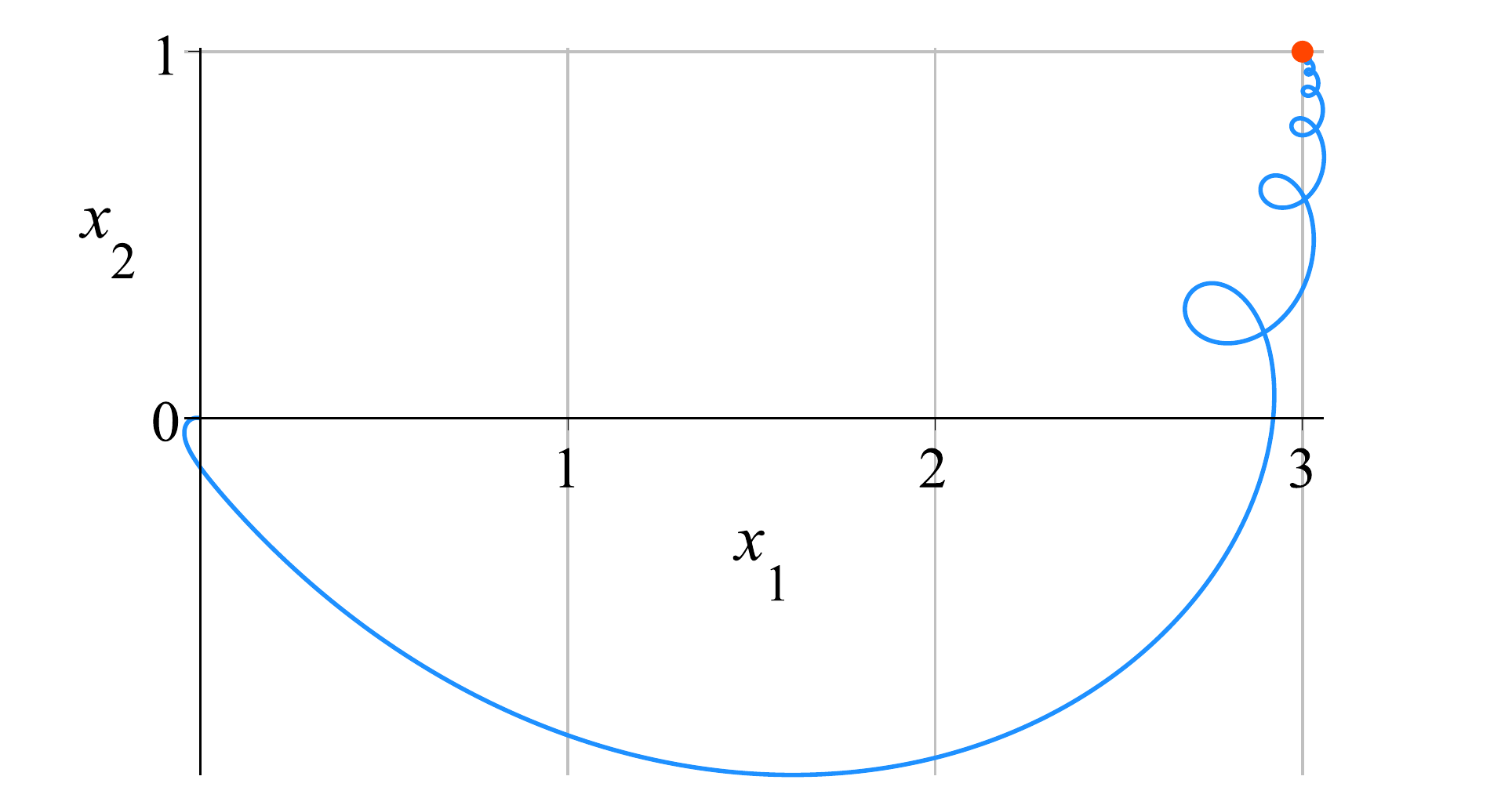}\\
\includegraphics[width=1\linewidth]{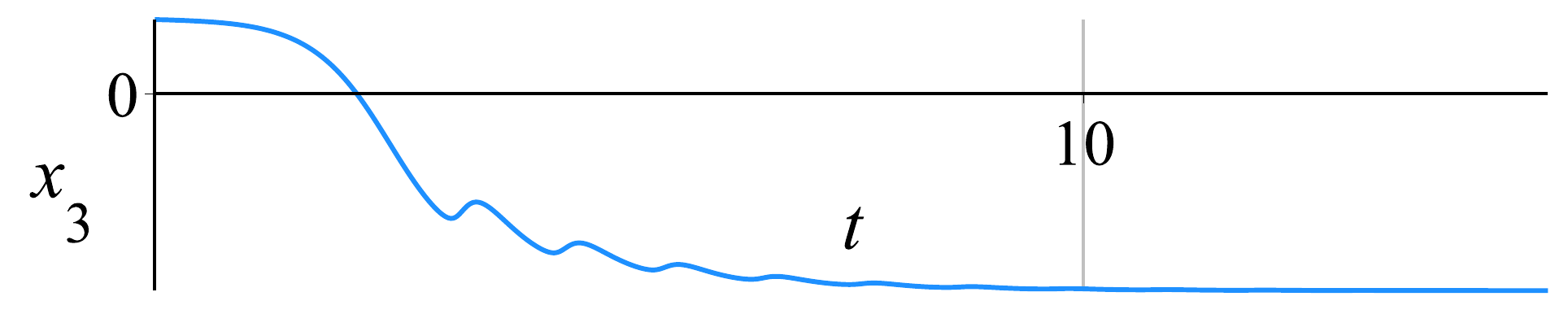}\\
b)
\end{center}
\end{minipage}
\begin{minipage}{0.32\linewidth}
\begin{center}
\includegraphics[width=1\linewidth]{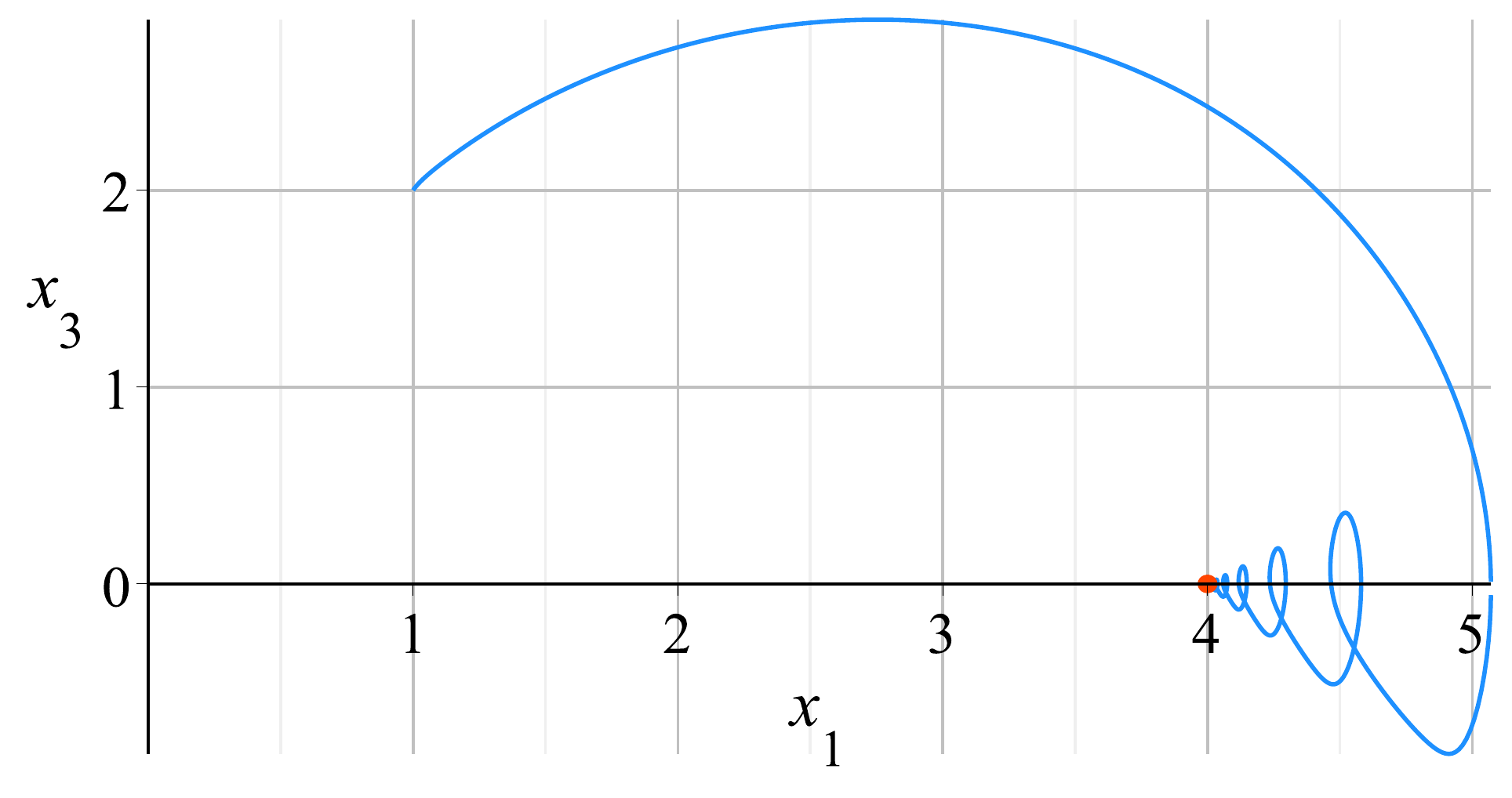}\\
\includegraphics[width=1\linewidth]{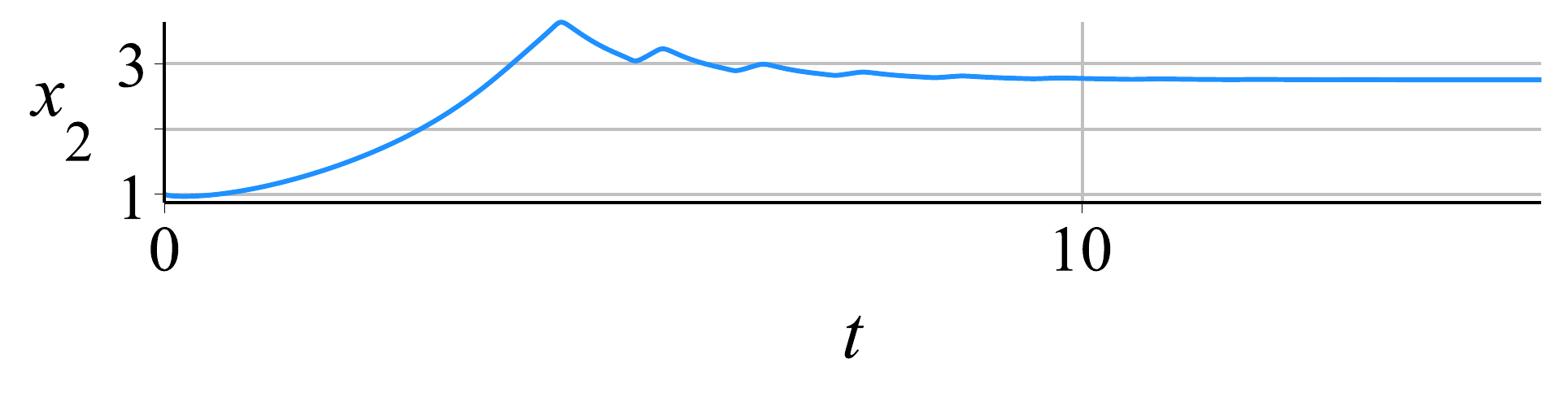}\\
c)
\end{center}
\end{minipage}
\caption{Projections of the trajectories of system~\eqref{ex_bro} on the  $(x_1,x_2)$-plane (top) and the graph of $x_3(t)$ (bottom) with controls~\eqref{cont_ex},\eqref{contA} (plot a)) and~\eqref{cont_ex},\eqref{contB} (plots b),c)). In the plots a),b), the cost function is given by~\eqref{Jx1x2}; $x(0)=(0,0,2)^\top$, $y^*=(3,1)^\top$.
In the plot c), the cost function is given by~\eqref{Jx1x3}; $x(0)=(1,1,2)^\top$, $y^*=(4,0)^\top$.
\label{bro_x1x2}
}
\end{figure*}
\begin{figure*}[tpt]
\includegraphics[width=0.33\linewidth]{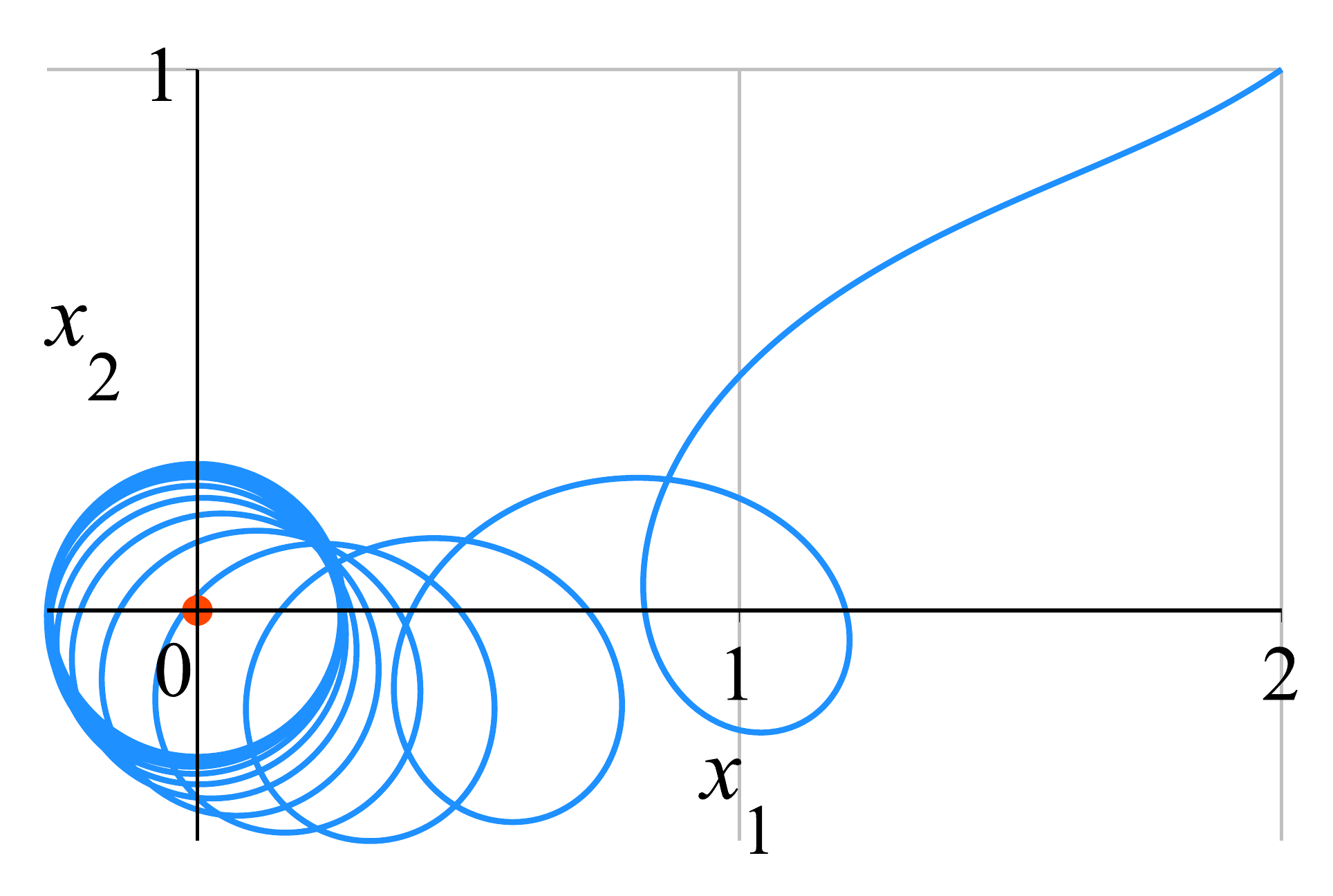}\hfill\includegraphics[width=0.33\linewidth]{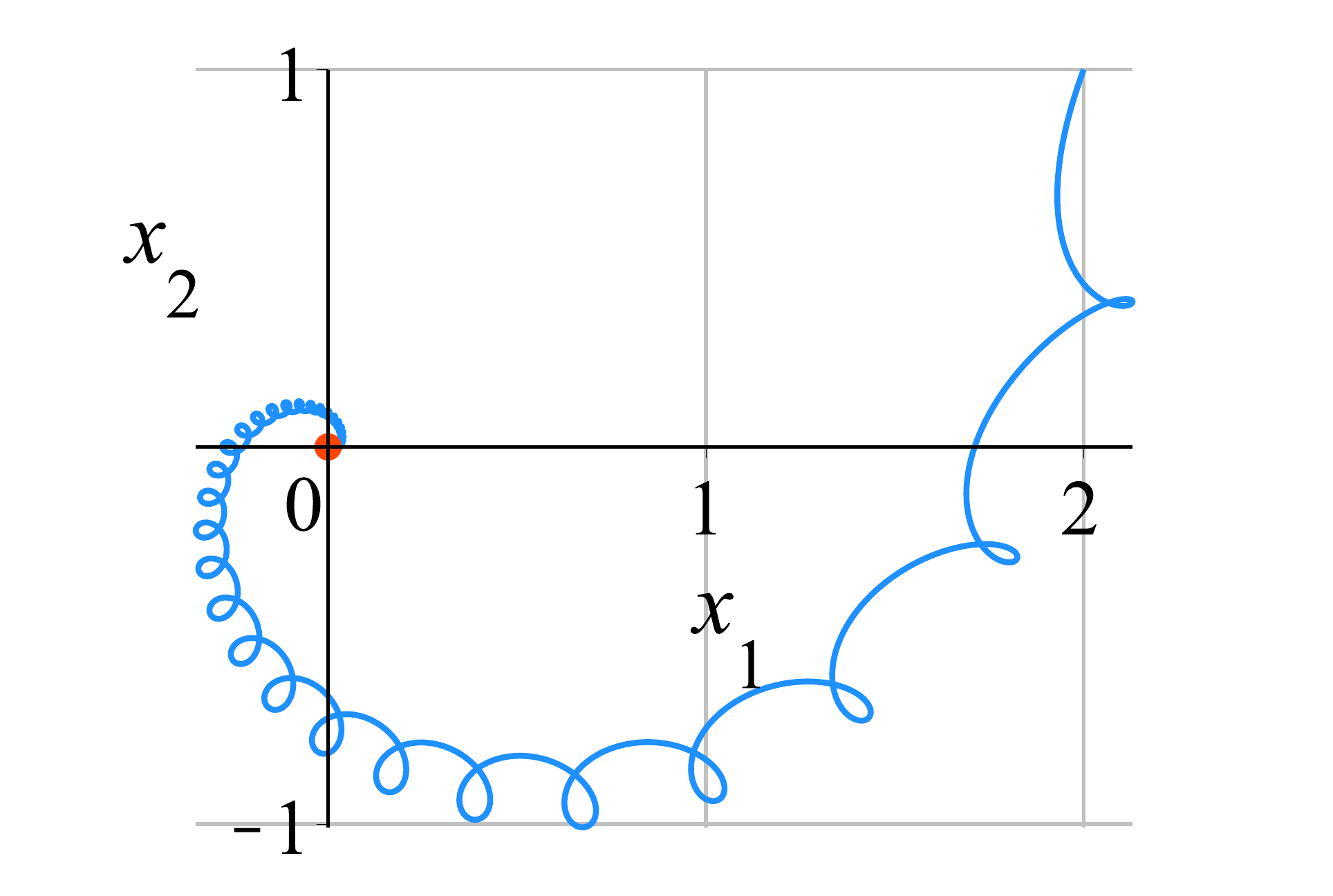} \hfill\includegraphics[width=0.33\linewidth]{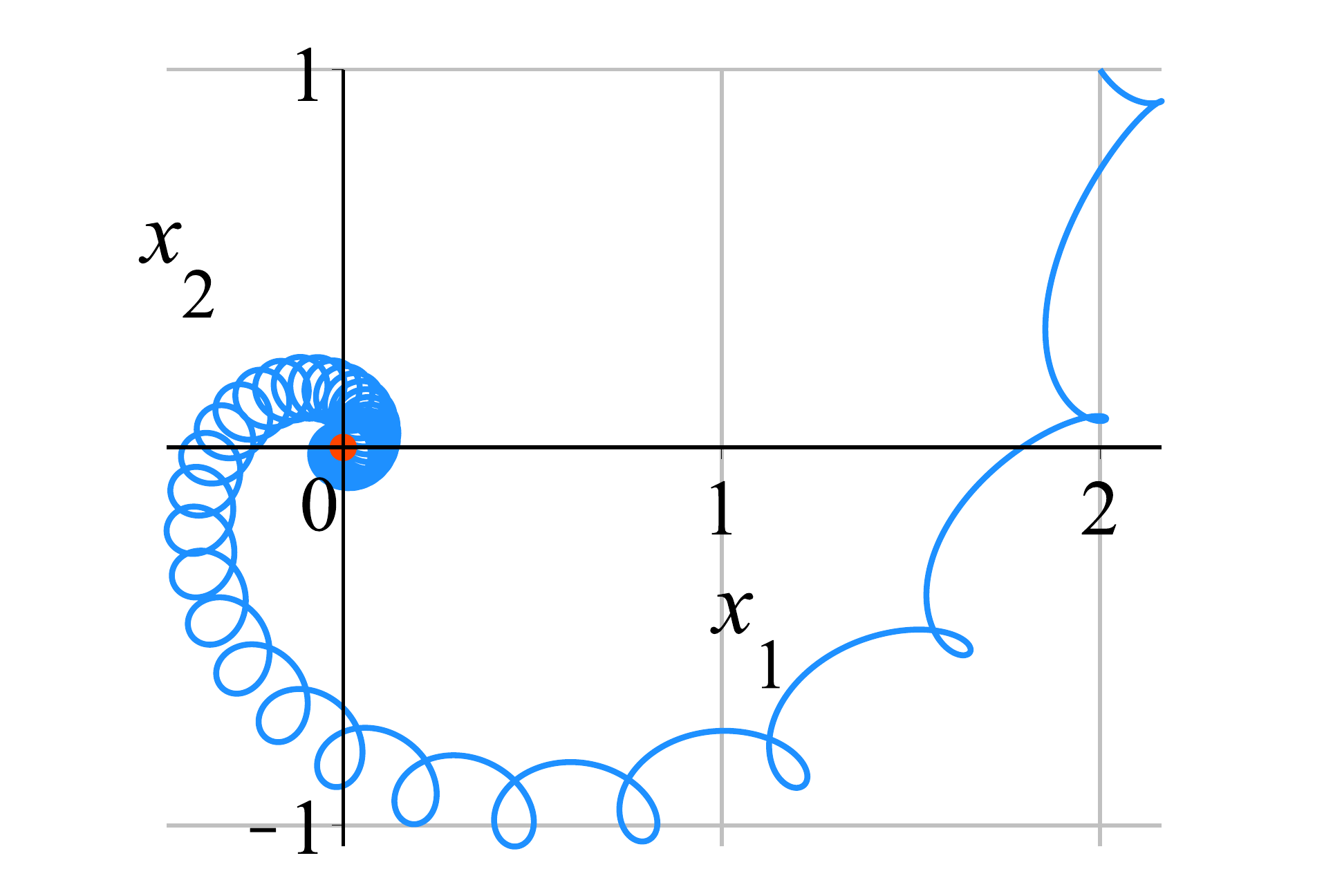}\\
\includegraphics[width=0.33\linewidth]{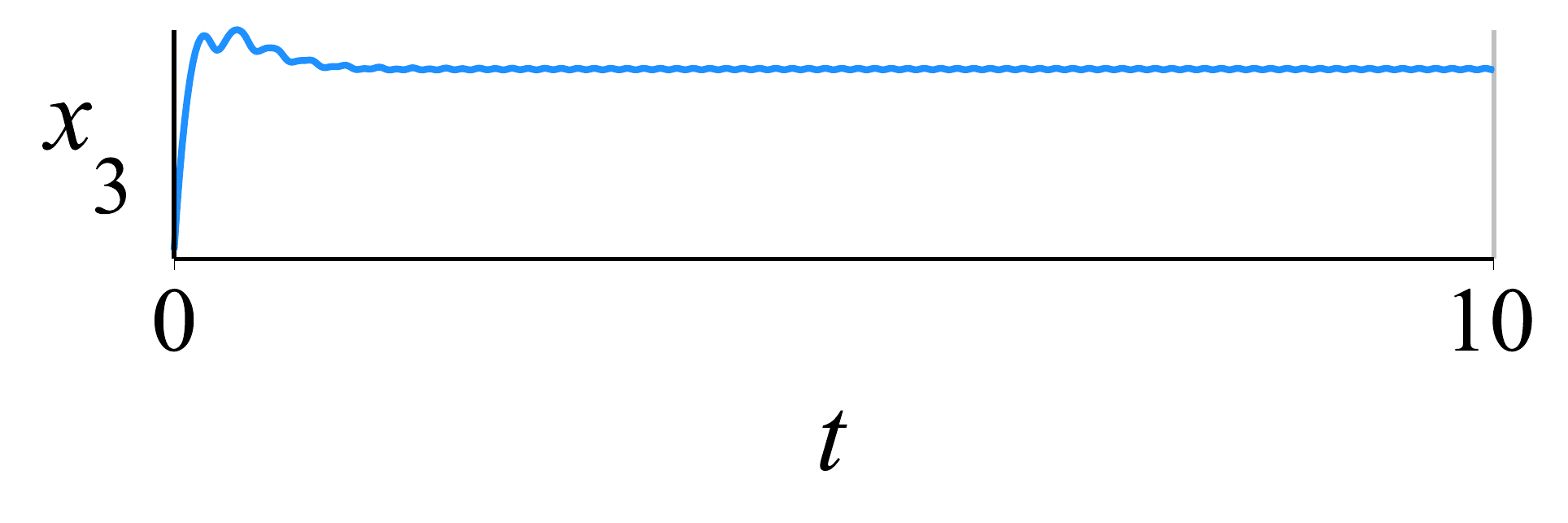}\hfill\includegraphics[width=0.33\linewidth]{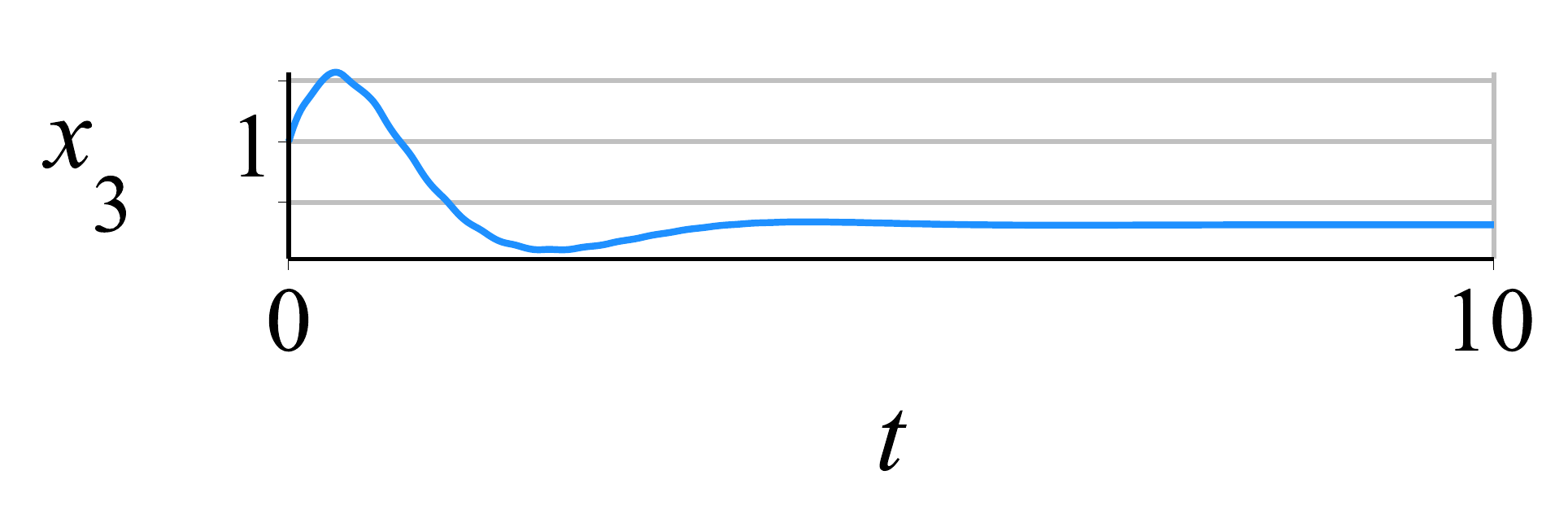} \hfill\includegraphics[width=0.33\linewidth]{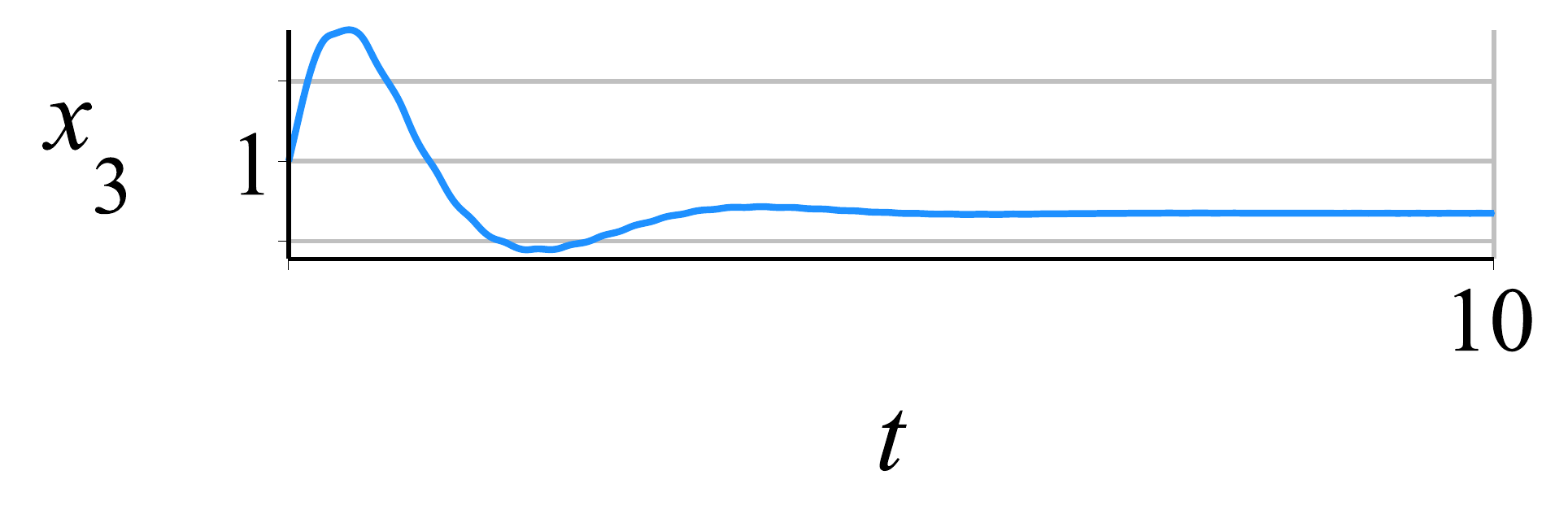}\\
a) \hfill b) \hfill c) \\
\caption{Projections of the trajectories of system~\eqref{body}  on the  $(x_1,x_2)$-plane (top) and  graph of $x_3(t)$ (bottom) with controls~\eqref{cont_ex},\eqref{contA} (plot a)) and ~\eqref{cont_ex},\eqref{contB} (plots b),c)), $J(x)=x_1^2+x_2^2$ (a),b)) and $J(x)=x_1^2+x_2^2+x_3^2$ (c)). Here $A_1=1,A_2=2,A_3=3$, $\varepsilon=0.25$, $x(0)=(2,1,1)^\top$, $y^*=(0 ,0)^\top$.\label{fig_body}}
\end{figure*}
\textbf{Sketch of the proof}. Computing the  time-derivative of  $J$ along the trajectories of the corresponding Lie bracket system for~\eqref{sys_es_1}, we get
$
\dot J(\bar y)=-\sum_{i,j=1}^{n_1}\gamma_i\Big(\frac{\partial J(\bar y)}{\partial \bar y_j}\tilde f_{ij}(\bar x)\Big)^2.
$
In general, $\dot J(\bar y)$ does not satisfy condition 1.2). However, it is easy to see that $\dot J(\bar y)=0$ if and only if
$
\nabla J(\bar y)F(\bar x)=0$, where $F(\bar x)=\left(
                                                      \begin{array}{ccc}
                                                        \tilde f_{11}(\bar x) & \dots &  \tilde f_{n_11}(\bar x) \\
                                                        \vdots & \ddots & \vdots \\
                                                         \tilde f_{1n_1}(\bar x) & \dots &  \tilde f_{n_1n_1}(\bar x) \\
                                                      \end{array}
                                                    \right).
$\\
Under the conditions of Theorem~\ref{thm_es2}, the matrix $F(\bar x)$ is nonsingular for all $\bar x$, which means  $\dot J(\bar y)=0$ if and only if $\bar y=y^*$. Then the practical asymptotic stability can be proved similar to Theorem~\ref{thm_gen_practical}. 
\section{Examples}
In this section, we consider several examples illustrating the obtained results and some possible extensions. In all the examples, we use extremum seeking controls $u=(u_1,u_2)^\top$ with
\begin{equation}\label{cont_ex}
\begin{aligned}
  u_1=\gamma_1\sqrt{\frac{\pi}{\varepsilon}}\Big(g_1(J(x))&\cos\frac{2\pi t}{\varepsilon} +g_{3}(J(x))\sin\frac{2\pi t}{\varepsilon}\Big),\\
  u_2=\gamma_2\sqrt{\frac{\pi}{\varepsilon}}\Big(g_2(J(x))&\sin\frac{2\pi t}{\varepsilon} -g_{4}(J(x))\cos\frac{2\pi t}{\varepsilon}\Big),
\end{aligned}
\end{equation}
where $\gamma_1,\gamma_2>0$, and the functions $g_i,g_{i+2}$ satisfy~\eqref{class}, $i=1,2$. We exploit two types of such functions:
\begin{equation}\label{contA}
g_i(z)=\sin z,\,g_{i+2}=\cos z;
\end{equation}
\begin{equation}\label{contB}
\begin{aligned}
g_i(z)&=\sqrt{\tfrac{1-e^{-z/4}}{1+e^{z/4}}}\sin(e^{z/4}+2\ln(e^{z/4}{-}1)),\\
g_{i+2}(z)&=\sqrt{\tfrac{1-e^{-z/4}}{1+e^{z/4}}}\cos(e^{z/4}+2\ln(e^{z/4}{-}1)),\,z>0,
\end{aligned}
\end{equation}
which were introduced in~(\cite{Sch14}) and~(\cite{GZE18}), respectively. Note that our reason for this is not to compare the performance of these control strategies, but just to illustrate different possibilities for control design.

\subsection{Partial stabilization of the Brockett integrator}
As the first example, we consider Problem~1 with the extremum seeking system described by the equations
\begin{equation}\label{ex_bro}
\begin{aligned}
  &\dot x_1=u_1,\ \dot x_2=u_2,\  \dot x_3=x_2u_1-x_1u_2,
\end{aligned}
\end{equation}
and the two cost functions:
\begin{equation}\label{Jx1x2}
  J_1(x_1,x_2)=(x_1-3)^2+(x_2-1)^2,
\end{equation}
\begin{equation}\label{Jx1x3}
  J_2(x_1,x_3)=(x_1-4)^2+x_3^2.
\end{equation}
For the  cost function $J_1(x_1,x_2)$, one can easily see that the assumptions of Theorem~\ref{thm_es2} are satisfied since the vector fields $\tilde f_1=(1,0)^\top$ and $\tilde f_2=(0,1)^\top$ are linearly independent in $\mathbb R^2$.
For $J_2(x_1,x_3)$,
$\tilde f_1=(1,x_2)^\top$ and $\tilde f_2=(0,-x_1)^\top$ are linearly independent if $x_1\ne 0$ which can be achieved if $x_1(0)x_1^*> 0$ and if $\varepsilon$ is small enough.
Note that the boundedness of the vector fields of~\eqref{ex_bro} holds only for controls~\eqref{cont_ex},\eqref{contB}, since in this case it can be proved that $x_1(t),x_2(t)$ belongs to a compact set for all $t\ge0$.
\\
Fig.~\ref{bro_x1x2},a) illustrates the behavior of trajectories of system~\eqref{ex_bro} with the cost function~\eqref{Jx1x2} and controls~\eqref{cont_ex},\eqref{contA}, $\varepsilon=0.75$, $\gamma_1=\gamma_2=2$. In this case, we observe  the practical asymptotic stability property. We expect that the use of controls~\eqref{cont_ex},\eqref{contB} yields the classical asymptotic stability result, similarly to the one obtained in~\cite{GZE18}. This property is illustrated in Fig.~\ref{bro_x1x2},b). For  the cost function~\eqref{Jx1x3}, the behavior of trajectories of system~\eqref{ex_bro} with controls~\eqref{cont_ex},\eqref{contB} is shown in Fig.~\ref{bro_x1x2},c).
\subsection{Partial stabilization of a rotating rigid body}
As another example,
consider the Euler  equations describing the rotational motion of a rigid body:
\begin{equation}\label{body}
\begin{aligned}
&\dot x_1=\tfrac{A_3-A_2}{A_1}x_2x_3+u_1\,,\dot x_2=\tfrac{A_1-A_3}{A_2}x_1x_3+u_2,\\
&\dot x_3=\tfrac{A_2-A_1}{A_3}x_1x_2.\\
\end{aligned}
\end{equation}
Here $x_1,x_2,x_3$ represent the principal components of the angular velocity vector,  $A_1,A_2,A_3>0$ are the main central moments of inertia, and $u_1,u_2$ are the control torques.
Our goal is to stabilize system~\eqref{body} along the $x_3$-axis, i.e. to $x_1^*=x_2^*=0$, assuming   that  the  cost function is $J(x)=x_1^2+x_2^2$. As in the previous example, we use  controls~\eqref{cont_ex},~\eqref{contA}, and~\eqref{contB}. Then the Lie bracket system for~\eqref{body} takes the form
\begin{align}
&\dot {\bar x}_1=\tfrac{A_3-A_2}{A_1}{\bar x}_2{\bar x}_3-2{\bar x}_1,\,\dot {\bar x}_2=\tfrac{A_1-A_3}{A_2}{\bar x}_1{\bar x}_3-2{\bar x}_2,\nonumber \\
&\dot {\bar x}_3=\tfrac{A_2-A_1}{A_3}{\bar x}_1{\bar x}_2. \label{body_lie}
\end{align}
Using the Lyapunov function $V({\bar x})=A_1{\bar x}_1^2+A_2{\bar x}_2^2+A_3{\bar x}_3^2$, one can show that $\dot V(\bar x)=-4(A_1x_1^2+A_2x_2^2)$. Note that in this case  condition 1.1) of Theorem~\ref{thm_gen_practical} is not satisfied; however, using Corollary~1 we can prove the practical asymptotic attractivity. Furthermore, if $\max\{A_1,A_2\}<A_3$ (or $\min\{A_1,A_2\}>A_3$), then the conditions of Theorem~\ref{thm_gen_practical} can be ensured with $V(x)=\tfrac{A_1}{A_3-A_2}x_1^2+\tfrac{A_2}{A_3-A_1}x_2^2$ (or $V(x)=\tfrac{A_1}{A_2-A_2}x_1^2+\tfrac{A_2}{A_1-A_3}x_2^2$) (see Fig.~\ref{fig_body},a) and b)).\\
The proposed techniques for generating partially stabilizing gradient-free controllers can also be used in related problems, e.g., for partial output stabilization of control systems. In particular, assume  that in the considered example only the measurements of
$J(x)=x_1^2+x_2^2+x_3^2$ are available. Then Corollary~1 implies that the controls~\eqref{cont_ex}, \eqref{contB} still can be used for steering system~\eqref{body} to a neighborhood of the  set $\{x{\in}\mathbb R^3:x_1{=}x_2{=}0\}$ (see Fig.~\ref{fig_body},c)).
\section{Conclusions}

In this paper, we have addressed the problem of extremum seeking with  respect to a part of variables. To obtain practical partial asymptotic stability conditions, we have extended the Lie bracket approximation approach and the methods proposed in~\cite{GZE18} to control-affine systems, whose averaged system has only a partially asymptotically stable equilibrium.
 The obtained results have been exploited for the design of extremum seeking controllers. Besides, we have illustrated  applications of the proposed techniques to partial output stabilization on the rotating rigid body example. In future work, we expect to derive classical (instead of practical) partial asymptotic stability conditions and relax assumptions on the Lyapunov function and the cost. Furthermore, we expect that the proposed approach will be of particular use for synchronization tasks.
\appendix
\section{Proof of Theorem~\ref{thm_gen_practical}}
Without loss of generality, assume $t_0=0$.
\\
For any $\delta\in\Big(0,\alpha_2^{-1}\big(\alpha_1({\rm dist}(y^*,\partial D_1))\big)\Big)$, let $c_\delta=\alpha_2(\delta)$, $D_0=\{(y,z)\in\mathbb R^n:\|y-y^*\|\le \delta,z\in D_2\}$. Then
$ D_0\subseteq D'=\{x: z\in D_2, V(x)\le c_\delta\}  \subset D.
$
From Assumption~1, we  define
\begin{equation}\label{as1}
  \begin{aligned}
  &M_0=\sup_{x\in D'}\|f_0(x)\|,\, M_1=\sup_{x\in D',1\le i\le m}\|f_i(x)\|\\
  &M_{2}=\sup_{x\in D',0\le i,j\le m}  \|L_{f_j}f_i(x)\|,\\
  &M_{3}=\sup_{\underset{0\le l\le m}{x\in D',1\le i,j\le m}}   \| L_{f_l} L_{f_j}f_i(x)\|.
\end{aligned}
\end{equation}
For any $\rho>0$, take $\delta' \in\Big(0,\alpha_2^{-1}\big(\alpha_1(\rho)\big)\Big)$ and put $\rho'=\alpha_1^{-1}\big(\alpha_2(\delta')\big)$,
$$d=\min\big\{\rho-\rho',{\rm dist}(y^*,\partial  D_1)-\alpha_1^{-1}(c_\delta)\big\}>0.$$
 By the conditions of Theorem~\ref{thm_gen_practical}, if $z(0)=z^0\in D_2$ then $z(t)\in D_2$ for all $t\ge 0$. Thus, to ensure that the solutions $x(t)$ with initial conditions $x(0)=x^0\in D'$ are well-defined in $D$ for $t\in[0,\varepsilon]$, it suffices to define $\varepsilon_0$ as the positive root of the equation $M_0\varepsilon+M_1W\sqrt\epsilon=d$. Then, for each $\varepsilon \in(0,\varepsilon_0)$, $x^0\in D'$, and for all $t\in[0,\varepsilon]$,
$$
\|y(t){-}y^*\|{\le }tM_0+\frac{tM_1W}{\sqrt\varepsilon}{+}\delta' {<}d{+}\delta'{<} {\rm dist}(y^*,\partial D_1).
$$
The above choice of $\rho',d$ implies the following properties:
\begin{equation}\label{rho}
\begin{aligned}
&  V(x^0){\le }\alpha_2(\delta'){\Rightarrow }\|y^0{-}y^*\|{<}\rho'{\Rightarrow} \|y(t){-}y^*\|{<}\rho,t\in[0,\varepsilon].\\
\end{aligned}
\end{equation}
To investigate the behavior of $V(x)$ along the trajectories of system~\eqref{sysA}, consider the Volterra series expansion of the solution $x(t)$ of system~\eqref{sysA} with an arbitrary initial condition $x(0)=x^0$ from $D'$ on the interval $t\in[0,\varepsilon]$:
\begin{equation}\label{volt1}
\begin{aligned}
  x(t)&=x^0+tf_0(x^0)+\frac{1}{\sqrt\varepsilon}\sum_{i=1}^mf_i(x^0)\int_0^tw_i\Big(\frac{\tau}{\varepsilon}\Big)d\tau\\
  &+\frac{1}{\varepsilon}\sum_{i<j}[f_i,f_j](x^0)\int_0^t\int_0^\tau\Big(w_j\Big(\frac{\tau}{\varepsilon}\Big)w_i\Big(\frac{s}{\varepsilon}\Big)\\
  &-w_i\Big(\frac{\tau}{\varepsilon}\Big)w_j\Big(\frac{s}{\varepsilon}\Big)\Big)dsd\tau+R(t),
\end{aligned}
\end{equation}
where
$$
\begin{aligned}
&R(t)=\int_0^t\int_0^\tau  L_{f_0}f_0(x(s))dsd\tau\\
&+ \frac{1}{\sqrt\epsilon}\sum_{i=1}^m\int_0^t\int_0^\tau\Big( L_{f_i}f_0(x(s))w_i\Big(\frac{s}{\varepsilon}\Big)+  L_{f_0}f_i(x(s))\\
&\times w_i\Big(\frac{\tau}{\varepsilon}\Big)\Big)dsd\tau+\frac{1}{\varepsilon}\sum_{i=1}^{m}L_{f_i}f_i(x^0)\int_0^t\int_0^\tau w_i\Big(\frac{\tau}{\varepsilon}\Big)\\
&\times w_i\Big(\frac{s}{\varepsilon}\Big)dsd\tau+\frac{1}{\varepsilon}\sum_{i,j=1}^{m}\int_0^t\int_0^\tau\int_0^s L_{f_0} L_{f_j}f_i(x(p))w_i\Big(\frac{\tau}{\varepsilon}\Big)\\
&\times w_j\Big(\frac{s}{\varepsilon}\Big)dpdsd\tau+\frac{1}{\varepsilon^{3/2}}\sum_{i,j,l=1}^{m}\int_0^t\int_0^\tau\int_0^s L_{f_l} L_{f_j}f_i(x(p))\\
&\times w_i\Big(\frac{\tau}{\varepsilon}\Big)w_j\Big(\frac{s}{\varepsilon}\Big)w_l\Big(\frac{p}{\varepsilon}\Big)dpdsd\tau.
\end{aligned}
$$
In particular, for $t=\varepsilon$,  representation~\eqref{volt1} takes the form
\begin{equation}\label{volt2}
  x(\varepsilon)=x^0+\varepsilon\Big( f_0(x^0)+\sum_{i<j}[f_i,f_j](x^0)\nu_{ij}\Big)+R(\varepsilon),
\end{equation}
and from~\eqref{as1} the remainder can be estimated as
$$
\|R(\varepsilon)\|{\le } 
\varepsilon^{3/2}\Big(M_2+\tfrac{W^2m^2M_3}{6}\Big)\Big({\sqrt\varepsilon}+Wm\Big){=}\sigma\varepsilon^{3/2},
$$
where $\sigma =\Big(M_2+\tfrac{W^2m^2M_3}{6}\Big)\Big({\sqrt\varepsilon}+Wm\Big)$ is monotone with respect to $\varepsilon$. Next, we apply Taylor's formula to  $V(x(\varepsilon))$:
$$
\begin{aligned}
V&(x(\varepsilon))=V(x^0)+\big(\nabla V(x^0),x(\varepsilon)-x^0\big)\\
&+\frac{1}{2}\sum_{i,j=1}^m\frac{\partial^2V(x)}{\partial x_i\partial x_j}\Big|_{x=x^0+\theta(x(\varepsilon)-x^0}(x_i(\varepsilon)-x_i^0)(x_j(\varepsilon)-x_j^0),
\end{aligned}
$$
with some $\theta\in(0,1)$. Let
$
\mu_1=\sup_{x\in D'}\big\|\nabla V(x)\big\|,
$
$
\mu_2=2\sup_{x\in D'}\Big\|\frac{\partial^2V(x)}{\partial x^2}\Big\|\big(M_0+M_2\sum_{i<j}\nu_{ij}+\sqrt\varepsilon\sigma\big)^2.
$
Then, from~\eqref{volt2} and~\eqref{as1}, we conclude that
$$
\begin{aligned}
V(x(\varepsilon))\le V(x^0)+\varepsilon L_{\bar f}V(x^0)+\varepsilon^{3/2}\sigma\mu_1+\varepsilon^2\mu_2.\\
\end{aligned}
$$
Recall that $ L_{\bar f}V(x)\le -\alpha_3(\|y-y^*\|)$ in $D$.  Thus, if $\|y^0 -y^*\|\ge\rho'$  then
 $$
\begin{aligned}
V(x(\varepsilon))\le V(x^0)&-\varepsilon\alpha_3(\rho')+\varepsilon^{3/2}\sigma\mu_1+\varepsilon^2\mu_2.
\end{aligned}
$$
Let $\lambda\in(0,\alpha_3(\rho'))$ and let $\varepsilon_1$ be the smallest positive root of the equation
$$
\sqrt\varepsilon\sigma\mu_1+\varepsilon^2\mu_2=\alpha_3(\rho')-\lambda.
$$
Then
\begin{equation}\label{Vdecay}
\begin{aligned}
V(x(\varepsilon))&\le V(x^0)-\varepsilon\lambda<V(x^0),
\end{aligned}
\end{equation}
provided that $\|y^0 -y^*\|\ge\rho'$.
The last inequality shows that $x(\varepsilon)\in D'$, and the solutions $x(t)$  of system~\eqref{sysA} with the initial conditions $x(0)=x^0\in D_0\subset  D'$ are well-defined in $D$ for $t\in[0,2\varepsilon]$.
Furthermore, we conclude that there exists an $N\in\mathbb N\cup \{0\}$ such that
\begin{equation}\label{yy}
\begin{aligned}
  & \|y(j\varepsilon) -y^*\|\ge\rho'\text{ for all }j=0,\dots,N-1,\\
 \text{and } & \|y(N\varepsilon) -y^*\|\le\rho'.
\end{aligned}
\end{equation}
Indeed, assume $\|y(j\varepsilon) -y^*\|\ge\rho'$ for all $j\in\mathbb N\cup\{0\}$. Then repeating inequality~\eqref{Vdecay}, we get
$
V(x(N\varepsilon))\le   V(x^0)-N\varepsilon\lambda.
$
With an increase of $N$, the right-hand side of the above inequality becomes negative which contradicts  $V(x(N\varepsilon)){\ge}0$. Thus, there exists an $N\in\mathbb N\cup \{0\}$ such that~\eqref{yy} holds.
\\
Estimate~\eqref{rho} implies that $\|y\big((N+1)\varepsilon\big)-y^*\|\le \rho$. If  $\|y\big((N+1)\varepsilon\big)-y^*\|\ge\rho'$, we apply~\eqref{Vdecay} again and obtain $$V\big((N+2)\varepsilon\big)<V\big((N+1)\varepsilon\big).$$
Otherwise we have $\|y\big((N+2)\varepsilon\big)-y^*\|\le \rho$ and repeat the procedure. Taking  $\bar\varepsilon=\min\{\varepsilon_0,\varepsilon_1\}$, we conclude that, for any $\varepsilon\in(0,\bar\varepsilon)$, the solutions of system~\eqref{sysA} satisfy the following property:
$$
\begin{aligned}
\text{if }&\|y(0)-y^*\|\le \delta\text{ and }z(0)\in D_2\text{ then there exists}\\
&\text{a }t_1>0\text{ such that }\|y(t)-y^*\|\le \rho \text{ for all }t\ge t_1.
\end{aligned}
$$
Since $\rho$ is assumed to be an arbitrary positive number, the practical $y$-attractivity has been proved. To prove the practical $y$-stability property, for any $\rho>0$ we take the $\delta'$ defined as before. Then, for any $y^0\in B_{\delta'}(y^*)\subset B_\rho(y^*)$ and $z^0\in D_2$, $V(x^0)\le\alpha_2(\delta')$. Summarizing~\eqref{rho},\eqref{Vdecay} and the previous argumentation, we  conclude with the stability property.
\end{document}